\documentclass[10pt]{amsart}

\usepackage{amssymb}
\usepackage{amsfonts}
\usepackage{amscd}
\usepackage{graphicx}
\usepackage{xypic}
\usepackage{psfrag}

\newcommand{\Gr}{\operatorname{Gr}}

\newcommand{\Hom}{\operatorname{Hom}}

\newcommand{\rank}{\operatorname{rank}}

\newcommand{\codim}{\operatorname{codim}}
\newcommand{\depth}{\operatorname{depth}}

\newcommand{\bull}{{\sssize \bullet}}
\newcommand{\scap}{{\,\ssize\cap\,}}
\renewcommand{\O}{\Omega}
\newcommand{\OO}{\boldsymbol{\Omega}}

\newcommand{\attach}[3]{(#1+#2,#3)}
\newcommand{\picA}[1]{\includegraphics{#1}}

\newcommand{\picC}[1]{\includegraphics[scale=0.50]{#1}}

\theoremstyle{plain}
\newtheorem*{main}{Main Theorem}
\newtheorem{lemma}{Lemma}
\newtheorem{cor}{Corollary}
\newtheorem*{conj}{Conjecture}
\newtheorem*{C1}{C1}
\newtheorem*{C2}{C2}
\newtheorem*{C3}{C3}

\numberwithin{equation}{section}

\newcommand{\refsec}[1]{Section~\ref{#1}}
\newcommand{\refeqn}[1]{(\ref{#1})}

\newcommand{\reflemma}[1]{Lemma~\ref{#1}}
\newcommand{\refcor}[1]{Corollary~\ref{#1}}

\title[Formulas for quiver varieties]{Chern class formulas for quiver
  varieties}
\date{\today}
\author{Anders Skovsted Buch}
\address{Department of Mathematics\\
  University of Chicago\\
  Chicago, IL 60637}
\email{abuch@math.uchicago.edu} 
\author{William Fulton}
\address{Department of Mathematics\\
  University of Chicago\\
  Chicago, IL 60637} 
\email{fulton@math.uchicago.edu}
\thanks{The research of the second author was supported by an Erlander
  Professorship in Sweden and the National Science Foundation}

\begin{document}

\psfrag{e-r}{$e - r$}
\psfrag{f-r}{$f - r$}

\psfrag{R}{$R$}

\psfrag{gam}{$\gamma$}
\psfrag{sig}{$\sigma$}
\psfrag{tau}{$\tau$}

\psfrag{R=}{$R =$}
\psfrag{Ti=}{$T_i =$}
\psfrag{Pi}{$P_i$}
\psfrag{Qi-1}{$Q_{i-1}$}
\psfrag{Ti-1i}{$T_{i-1,i}$}

\psfrag{R01}{$R_{01}$}
\psfrag{R12}{$R_{12}$}
\psfrag{R23}{$R_{23}$}
\psfrag{Rn-1n}{$R_{n-1,n}$}
\psfrag{R02}{$R_{02}$}
\psfrag{R13}{$R_{13}$}
\psfrag{Rn-2n}{$R_{n-2,n}$}

\psfrag{e0-r1}{$e_0 - r_1$}
\psfrag{e1-r1}{$e_1 - r_1$}
\psfrag{e1-r2}{$e_1 - r_2$}
\psfrag{e2-r2}{$e_2 - r_2$}
\psfrag{e2-r3}{$e_2 - r_3$}
\psfrag{e3-r3}{$e_3 - r_3$}
\psfrag{en-1-rn}{$e_{n-1} - r_n$}
\psfrag{en-rn}{$e_n - r_n$}
\psfrag{r1}{$r_1$}
\psfrag{r2}{$r_2$}
\psfrag{r3}{$r_3$}
\psfrag{rn-1}{$r_{n-1}$}
\psfrag{rn}{$r_n$}

\psfrag{r00}{$r_{00}$}
\psfrag{r11}{$r_{11}$}
\psfrag{r22}{$r_{22}$}
\psfrag{r33}{$r_{33}$}
\psfrag{r44}{$r_{44}$}
\psfrag{r01}{$r_{01}$}
\psfrag{r12}{$r_{12}$}
\psfrag{r23}{$r_{23}$}
\psfrag{r34}{$r_{34}$}
\psfrag{r02}{$r_{02}$}
\psfrag{r13}{$r_{13}$}
\psfrag{r24}{$r_{24}$}
\psfrag{r03}{$r_{03}$}
\psfrag{r14}{$r_{14}$}
\psfrag{r04}{$r_{04}$}

\psfrag{rij}{$r_{ij}$}
\psfrag{ri-1j}{$r_{i-1,j}$}
\psfrag{rij+1}{$r_{i,j+1}$}
\psfrag{rij-1}{$r_{i,j-1}$}
\psfrag{ri+1j}{$r_{i+1,j}$}

\maketitle

\section{Introduction}

Our goal in this paper is to prove a formula for the general
degeneracy locus $\O_r$ associated to an oriented quiver of type
$A_n$.  If we are given a sequence of vector bundles and vector
bundle maps
\[
E_0 \xrightarrow{\phi_1} E_1 \xrightarrow{\phi_2} E_2 \to \dots \to
E_{n-1} \xrightarrow{\phi_n} E_n
\]
on an algebraic variety $X$, and a
collection $r = (r_{ij})_{0 \leq i < j \leq n}$ of non-negative
integers, there is a degeneracy locus $\O_r = \O_r(E_{\bull}) =
\O_r(E_{\bull},\phi_{\bull})$ defined by
\begin{equation} \label{1.1}
\O_r = \{x \in X \mid \rank\bigl(E_i(x) \to E_j(x)\bigr) \leq
r_{ij} ~\forall~ i < j \} \,.
\end{equation}
This is a closed subscheme of $X$; locally, where the bundles are
trivial, this is defined by vanishing of the minors of size $r_{ij}+1$
in the product of matrices giving the map $\phi_j \circ \cdots \circ
\phi_{i+1}$ from $E_i$ to $E_j$, for all $i < j$.

Not all rank conditions give reasonable loci.  Those that do---and the
only ones we will consider---are characterized by the conditions
\begin{equation}
\label{eqn_rankcond}
\begin{aligned}
r_{ij} \leq r_{i,j-1} \quad \text{and} \quad r_{ij} \leq r_{i+1,j} 
\quad \text{for all } i < j \text{, and} \\
r_{i+1,j-1} - r_{i,j-1} - r_{i+1,j} + r_{ij} \geq 0 \quad \text{for all }
i < j - 1,
\end{aligned}
\end{equation}
where we set $r_{ii} = \rank(E_i)$.  In fact, rank conditions
satisfying \refeqn{eqn_rankcond} are the only conditions that can
actually \emph{occur}, i.e.\ for which one can have equality in
\refeqn{1.1}.  When the maps are sufficiently generic, each such
$\O_r$ is irreducible, of codimension
\begin{equation}\label{1.3}
d(r) = \sum_{i < j} (r_{i,j-1} - r_{ij})(r_{i+1,j} - r_{ij}) \,.
\end{equation}

When $n = 1$, the formula for $\O_r$ is the well-known
Giambelli-Thom-Porteous formula, which we recall in order to introduce
some notation.  For a map $\phi : E \to F$ of vector bundles of ranks
$e$ and $f$, and a non-negative integer $r \leq \min(e,f)$, $\O_r$ is
the locus where $\phi$ has rank at most $r$.  The formula for $\O_r$
is the Schur polynomial
\[
s_{(e-r)^{f-r}}(F - E) \,,
\]
which is defined as follows.  Define cohomology classes $h_i$ by the
formula $\sum h_i = c(E^\vee)/c(F^\vee)$, where $ c(E^\vee) = 1 -
c_1(E) + c_2(E) - \cdots$ is the total Chern class, and the division
is carried out formally; in particular, $h_0 = 1$ and $h_i = 0$ for $i
< 0$.  For any sequence $\lambda = (\lambda_1, \dots , \lambda_p)$ of
non-negative integers, set
\[
s_\lambda(F - E) = \det(h_{\lambda_i+j-i})_{1 \leq i,j \leq p} \,.
\]
In the Giambelli-Thom-Porteous formula, $\lambda = (e-r)^{f-r}$
denotes the sequence $e-r$ repeated $f-r$ times.  In a Schur
determinant $s_\lambda(F - E)$, $\lambda$ will usually be a partition,
i.e.\ a weakly decreasing sequence, but later we will also need this
notation when $\lambda$ is not a partition.

Our general formula for the locus $\O_r$, when $r$ is any set of rank
conditions satisfying \refeqn{eqn_rankcond}, has the form
\[ \sum_\lambda c_\lambda(r) \, s_\lambda(E_\bull) \,, \]
where the sum is over sequences $\lambda = (\lambda(1), \lambda(2),
\dots, \lambda(n))$, with each $\lambda(i)$ a partition.  The class
$s_\lambda(E_\bull)$ is defined to be 
\[
s_\lambda(E_\bull) = s_{\lambda(1)}(E_1 - E_0) \cdot s_{\lambda(2)}(E_2 -
E_1) \cdot \ldots \cdot s_{\lambda(n)}(E_n - E_{n-1}) \,.
\]
The coefficients $c_\lambda(r)$ are certain integers for which we give
an inductive formula.

A second purpose of this paper is to introduce these integers
$c_\lambda(r)$, which we regard as generalized Littlewood-Richardson
coefficients.  We have a conjectured formula for $c_\lambda(r)$ as the
number of sequences $(T_1, T_2, \dots, T_n)$ of Young tableaux, with
$T_i$ of shape $\lambda(i)$, satisfying certain conditions.  This
formula has been proved when the number of bundles is at most four,
but it appears to be a difficult combinatorial problem to prove it in
general.

In \cite{fulton:universal} a special case of this situation was studied,
where the rank conditions are given by a permutation $w$.  For maps
\[ G_1 \to G_2 \to \dots \to G_m \to F_m \to F_{m-1} \to \dots \to F_1 \]
with $\rank(G_i) = \rank(F_i) = i$, and $w \in S_{m+1}$, let
\[ \O_w = \{ x \in X \mid \rank(G_q(x) \to F_p(x)) \leq r_w(p,q)
   ~\forall~ p,q \leq m \} \,,
\]
where $r_w(p,q) = \#\{ i \leq p \mid w(i) \leq q \}$.  These loci are
special cases of the loci $\O_r$ described in this paper.  The
formulas given here therefore specialize to the \emph{universal double
Schubert polynomials} $\mathfrak{S}_w(c_\bull(F_\bull);
c_\bull(G_\bull))$ for these loci.  Since these universal Schubert
polynomials specialize to quantum and double Schubert polynomials
(\cite{fomin.gelfand.ea:quantum}, \cite{kirillov.maeno:quantum},
\cite{fulton:flags}), we derive formulas for
these important polynomials.  These formulas appear to be new even for
the single Schubert polynomials $\mathfrak{S}_w(x)$.

Among the loci considered here are the \emph{varieties of complexes},
which are the loci $\O_r$ with $r_{ij} = 0$ for $j-i \geq 2$.  In this
case the formula for the coefficients $c_\lambda(r)$ is particularly
simple, and it agrees with our general conjectured formula.

In \refsec{sec_quiver} we discuss the loci $\O_r$ in more detail,
state the main theorem, and derive the main applications. This
includes a precise statement of what it means for a polynomial $\sum
c_\lambda(r)s_\lambda(E_\bull)$ to give a formula for a locus $\O_r$.
This statement implies the assertion that if $X$ is non-singular and
$\O_r$ has the expected codimension $d(r)$, then
\begin{equation}
\label{eqn_main}
[\O_r] = \sum_\lambda c_\lambda(r) \, s_\lambda(E_\bull) 
\end{equation}
in the Chow group $A^{d(r)}(X)$.  However, weaker assertions can be
made when $X$ is singular or the maps $\phi_i$ are less generic.  At
the end of \refsec{sec_proof} we sketch a generalization, which is
based on explicit resolutions of singularities of these loci.

The coefficients $c_\lambda(r)$ are determined by the geometry, if
this assertion is interpreted correctly.  We will see in
\refsec{sec_quiver} that $c_\lambda(r)$ depends only on the
differences $r_{i,j-1} - r_{ij}$ and $r_{i+1,j} - r_{ij}$.  This
allows the ranks of the bundles $E_i$ to be taken large compared to
the expected codimension $d(r)$; if the Chern classes of the bundles
are independent, the coefficients $c_\lambda(r)$ are then uniquely
determined by \refeqn{eqn_main}.

Much of the work in a project of this kind---discovering the shape of
the formula---is invisible in the final product, which has a short
proof (given in \refsec{sec_proof}).  In particular, it came as a pleasant
surprise to us that the polynomials for all the loci $\O_r$ can be
written as a linear combination of the polynomials
$s_\lambda(E_\bull)$.  We know of no reason for this other than the
proof of the explicit formula.  That the coefficients $c_\lambda(r)$
appear to be non-negative is even more surprising.

The conjectured formula for the coefficients is discussed in more
detail in the final \refsec{sec_coefs}; proofs of the combinatorial
assertions made there can be found in \cite{buch:combinatorics}.

We are particularly grateful to S.~Fomin, who provided an involution
on pairs of tableaux which gave us the strongest evidence for the
conjectured formula, and who has collaborated with us on the
combinatorial aspects of this problem.  Thanks also to M.~Haiman and
M.~Shimozono for their responses to combinatorial questions.
The Schubert package \cite{katz.strmme:schubert} was useful for
calculations.

\section{Quiver varieties; the theorem and applications}
\label{sec_quiver}

\subsection{The Main Theorem}
\label{ss_main}

Given vector bundles $E_0, \ldots, E_n$ on a variety $X$, let $H$ be
the direct sum of the bundles $\Hom(E_{i-1}, E_i)$, i.e.
\[
H = \Hom(E_0, E_1) \times_X \Hom(E_1, E_2) \times_X \dots 
    \times_X \Hom(E_{n-1}, E_n) \,.
\]
Writing $\widetilde E_i$ for the pullback of $E_i$ to $H$, we have a
universal or tautological sequence of bundle maps
\begin{equation}
\label{eqn_tautological}
\widetilde E_0 \xrightarrow{\Phi_1} \widetilde E_1
\xrightarrow{\Phi_2} \widetilde E_2 \to \dots \to \widetilde E_{n-1}
\xrightarrow{\Phi_n} \widetilde E_n
\end{equation}
on $H$.  For this universal case, it is a theorem of Lakshmibai and
Magyar \cite{lakshmibai.magyar:degeneracy} that for $r$ satisfying
\refeqn{eqn_rankcond}, the scheme $\widetilde \O_r = \O_r(\widetilde
E_\bull)$ for $(\ref{eqn_tautological})$ is reduced and irreducible,
of codimension $d(r)$.  Moreover, $\widetilde \O_r$ is a
Cohen-Macaulay variety if $X$ is Cohen-Macaulay.  (Earlier Abeasis,
del Fra, and Kraft \cite{abeasis.del-fra.ea:geometry} had shown, in
characteristic zero, that the reduced scheme $(\widetilde
\O_r)_{\text{red}}$ is Cohen-Macaulay.)  Note that, when the bundles
are trivial, $H$ is a Cartesian product of $X$ and a product $M$ of
spaces of matrices, and $\widetilde \O_r$ is a product of $X$ with the
corresponding locus in $M$; it is this locus in $M$ that is studied in
\cite{abeasis.del-fra:degenerations},
\cite{abeasis.del-fra.ea:geometry}, and
\cite{lakshmibai.magyar:degeneracy}.

The statement that ``the polynomial $P = \sum c_\lambda(r)
s_\lambda(E_\bull)$ is a formula for the locus $\O_r$'' has the usual
meaning in intersection theory (cf.\ \cite[\S14]{fulton:intersection},
\cite[App.~A]{fulton.pragacz:schubert}).  It implies that when $X$ is
non-singular and $\codim(\O_r, X) = d(r)$, then
\[ [\O_r] = \sum_\lambda c_\lambda(r) \, s_\lambda(E_\bull) \]
in the Chow group $A^{d(r)}X$, where $[\O_r]$ is the cycle defined
by the scheme $\O_r$.  For arbitrary $X$ and maps $\phi_i$,
there is a well defined cycle class $\OO_r$ in the Chow group
$A_{m - d(r)}(\O_r)$, where $m = \dim(X)$, whose image in $A_{m -
  d(r)}(X)$ is the class $\sum c_\lambda(r) s_\lambda(E_\bull) \scap
[X]$.  Whenever $\O_r$ has codimension $d(r)$ in $X$, $\OO_r$ is a
positive cycle supported on $\O_r$; if $X$ is Cohen-Macaulay, or more
generally if $\depth(\Omega_r, X) = d(r)$, this cycle is $[\O_r]$, but
if $X$ is not Cohen-Macaulay the coefficient of a component of $\O_r$
in $\OO_r$ may be smaller than the length of $\O_r$ at its generic
point.  These classes $\OO_r$ are compatible with the basic
constructions of intersection theory, exactly as in
\cite[Thm.~14.3]{fulton:intersection}.

In fact, to give maps $\phi_i : E_{i-1} \to E_i$ for all $i$
is the same as giving a section $s : X \to H$ of the bundle
$H$, and $\O_r = s^{-1}(\widetilde \O_r)$.  The general class $\OO_r$ is
constructed by intersecting $\widetilde \O_r \subset H$ with the (regular)
embedding $s : X \to H$, i.e.
\[ \OO_r = s^![\widetilde \O_r] \,, \]
where $s^! : A_*(\widetilde \O_r) \to A_*(\O_r)$ is the refined
intersection \cite[\S6]{fulton:intersection}.  As in
\cite[\S14]{fulton:intersection}, the general properties of these
classes
follow from this construction.  It therefore suffices to prove the
corresponding formula on $H$, i.e.\ that 
\[ [\widetilde \Omega_r] = 
   \sum_\lambda c_\lambda(r) \, s_\lambda(\widetilde E_\bull) \scap [H]
\]
in $A_{N - d(r)}(H)$, where $N = \dim(H) = \dim(X) + \sum_{i=1}^n
e_{i-1} e_i$, $e_i = \rank(E_i)$.

It is natural to arrange the rank conditions in a triangular array:
\setcounter{MaxMatrixCols}{20}
\[ \begin{matrix}
r_{00} && r_{11} && r_{22} && r_{33} && \cdots && r_{nn} \\
& r_{01} && r_{12} && r_{23} && \cdots && r_{n-1,n} \\
&& r_{02} && r_{13} && \cdots && r_{n-2,n} \\
&&& r_{03} && \cdots && r_{n-3,n} \\
&&&& \ddots \\
&&&&& r_{0n}
\end{matrix} \]
It is useful to replace each small triangle 
\[ \begin{matrix} e && f \\ & r \end{matrix} \]
occurring in this array by the rectangle of width $e-r$ and height
$f-r$.
\[ \picC{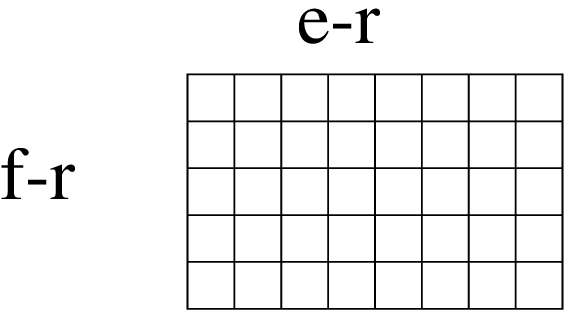} \]
We then have the rectangular array
\[ \begin{matrix}
R_{01} && R_{12} && R_{23} && \cdots && R_{n-1,n} \\
& R_{02} && R_{13} && \cdots && R_{n-2,n} \\
&& R_{03} && \cdots && R_{n-3,n} \\
&&& \ddots \\
&&&& R_{0n}
\end{matrix} \]
where $R_{ij}$ has width $r_{i,j-1} - r_{ij}$ and height $r_{i+1,j} -
r_{ij}$.  Note that the expected codimension $d(r)$ is the sum of the
areas of the rectangles.  The condition \refeqn{eqn_rankcond} says
that the rectangles get (weakly) shorter as one proceeds in a
southeasterly direction, and they get (weakly) narrower as one travels
southwest.  For example, the rank conditions given in the triangular
array
\[ \begin{matrix}
6 && 8 && 9 && 6 \\
& 5 && 6 && 6 \\
&& 4 && 3 \\
&&& 2
\end{matrix} \]
correspond to the rectangular array:
\[ \picC{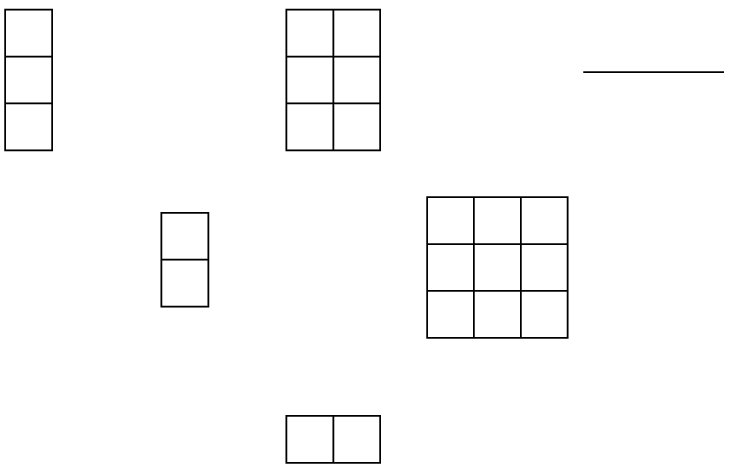} \]
Our formula depends on the rectangles in this array.  To be precise,
it depends on the integers $r_{i,j-1}-r_{ij}$ and $r_{i+1,j}-r_{ij}$
for all $i < j$; if a width $r_{i,j-1}-r_{ij}$ is zero, we need to
know the height $r_{i+1,j} - r_{ij}$, even though the rectangle
$R_{ij}$ is empty.  (The conjectured formula discussed later does not
have this defect.)  Each $R_{ij}$ is identified with the partition
$(r_{i,j-1} - r_{ij})^{r_{i+1,j} - r_{ij}}$ for which it is the Young
diagram.

At this point we need some notation.  If $R$ is a rectangle of width
$e$ and height $f$, and $\sigma$ and $\tau$ are partitions, with the
length $\ell(\sigma)$ of $\sigma$ at most $f$, then
$\attach{R}{\sigma}{\tau}$ denotes the sequence $(e + \sigma_1, e +
\sigma_2, \dots, e + \sigma_f, \tau_1, \tau_2, \dots)$; this is a
partition only if $e + \sigma_f \geq \tau_1$.  For a partition
$\lambda = (\lambda_1, \dots, \lambda_p)$, $|\lambda|$ denotes $\sum
\lambda_i$, which is the number of boxes in the Young diagram of
$\lambda$.  For partitions $\lambda$, $\sigma$, $\tau$ with $|\sigma|
+ |\tau| = |\lambda|$, $c^\lambda_{\sigma,\tau}$ denotes the
Littlewood-Richardson coefficient, which is the coefficient of the
Schur polynomial $s_\lambda$ in the expansion of $s_\sigma \cdot
s_\tau$ (see \cite{macdonald:symmetric*2}).  We set $e_i = r_{ii} =
\rank(E_i)$, $r_i = r_{i-1,i}$, and $R_i = R_{i-1,i}$, so $R_i$ has
height $e_i - r_i$ and width $e_{i-1} - r_i$.

If $I = (i_1, \dots, i_p)$ is a sequence of non-negative integers that
is not weakly decreasing, then $s_I(F-E)$ is either $0$ or it is $\pm
s_\lambda(F-E)$ for some unique partition $\lambda$ and unique
coefficient $\pm 1$.  This partition and coefficient can be found by
performing a sequence of moves of the type
\[ (j_1, \dots, j_p) \mapsto 
(j_1, \dots, j_{k-1}, j_{k+1} - 1, j_k + 1, j_{k+2}, \dots, j_p)
\]
if $j_{k+1} > j_k$; if one reaches a sequence $(j_1, \dots, j_p)$ with
some $j_{k+1} = j_k + 1$, then $s_I = 0$; otherwise one reaches a
partition $\lambda = (\lambda_1, \dots, \lambda_p)$ in $m$ steps, and
then $s_I = (-1)^m s_\lambda$.

We now give an algorithm for constructing finite formal sums $\sum
c_\lambda(r) S(\lambda)$, with $\lambda$ varying over $n$-tuples of
partitions $\lambda = (\lambda(1), \dots, \lambda(n))$.  The
polynomial for the degeneracy locus $\O_r$ will be obtained by
replacing each $S(\lambda)$ by $s_\lambda(E_\bull)$.  In the algorithm
we will meet symbols $S(I(1), \dots, I(n))$ where each $I(j)$ is a
sequence of integers.  For such symbols we imitate the above rule for
Schur polynomials to write $S(I(1), \dots, I(1))$ as either zero or
$\pm S(\lambda(1), \dots, \lambda(n))$ for unique partitions
$\lambda(1), \dots, \lambda(n)$.  If $s_{I(j)} = 0$ for any $j$, put
$S(I(1), \dots, I(n)) = 0$; otherwise write $s_{I(j)} = \epsilon_j
s_{\lambda(j)}$ for $1 \leq j \leq n$, with $\epsilon_j = \pm 1$, and
put $S(I(1), \dots, I(n)) = (\prod \epsilon_j) S(\lambda(1), \dots,
\lambda(n))$.

We construct the polynomial $\sum c_\lambda(r) S(\lambda)$ by
induction on $n$.  For $n = 1$ we have just one rectangle $R =
R_{01}$, and the polynomial is $S(R)$, which gives $s_R(E_1 - E_0)$.
Given the rectangular array for $r$, delete the top row.  This gives a
smaller array, for which we have a polynomial $\sum d_\mu S(\mu)$ by
induction, the sum over sequences $\mu = (\mu(1), \dots, \mu(n-1))$ of
partitions.  The polynomial $\sum c_\lambda(r) S(\lambda)$ is obtained
by replacing each $S(\mu)$ in $\sum d_\mu S(\mu)$ by
\[ \sum \left( \prod_{i=1}^{n-1} c^{\mu(i)}_{\sigma(i),\tau(i)} \right)
   S(I(1), \dots, I(n)) \,.
\]
Here the sum is over all sequences $(\sigma(1), \dots, \sigma(n-1))$
and $(\tau(1), \dots, \tau(n-1))$ of partitions, with $|\sigma(i)| +
|\tau(i)| = |\mu(i)|$, such that the length of $\sigma(i)$ is at most
the height of $R_i$, i.e.\ $\ell(\sigma(i)) \leq e_i - r_i$.
Define $I(i)$ to be $\attach{R_i}{\sigma(i)}{\tau(i-1)}$ for $1 \leq i
\leq n$, where $\tau(0)$ and $\sigma(n)$ are taken to be the empty
partition. 
One uses the rules just given to write each $S(I(1), \dots, I(n))$ as
$0$ or $\pm S(\lambda)$ for a unique $\lambda = (\lambda(1), \dots,
\lambda(n))$, thus arriving at a polynomial $\sum c_\lambda(r)
S(\lambda)$.

\begin{main}
The formula for $\Omega_r$ is $\sum c_\lambda(r) s_\lambda(E_\bull)$.
\end{main}

This theorem will be proved in the next section.  We first interpret
it in the case where the rectangular array has only two non-empty
rows, i.e.\ $R_{ij} = \emptyset$ if $j-i > 2$.  In this case the
inductive polynomial $\sum d_\mu S(\mu)$ is just $S(\mu)$, for $\mu =
(R_{02}, R_{13}, \dots, R_{n-2,n})$.  For a rectangle shape $R$, the
Littlewood-Richardson coefficient $c^R_{\sigma,\tau}$ vanishes unless
$\sigma$ and the $180^\circ$ rotation of $\tau$ {\em fit together\/}
to make $R$, in which case $c^R_{\sigma,\tau} = 1$.
\[ \picC{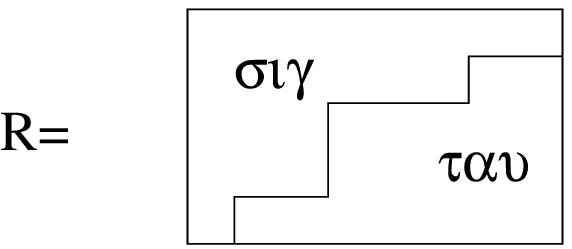} \]

\begin{cor}
\label{cor_complexes}
If $R_{ij}$ is empty for $j-i > 2$, then the formula for $\O_r$ is
$\sum s_\lambda(E_\bull)$, where the sum is over all $\lambda =
(\lambda(1), \dots, \lambda(n))$, with $\lambda(i) =
\attach{R_i}{\sigma(i)}{\tau(i-1)}$, such that $\sigma(i)$ and
$\tau(i)$ fit together to form $R_{i-1,i+1}$ for $1 \leq i \leq n-1$;
here $\sigma(n)$ and $\tau(0)$ are empty.
\end{cor}

Note that, by \refeqn{eqn_rankcond}, for any division of $R_{i-1,i+1}$
into $\sigma(i)$ and $\tau(i)$, $\sigma(i)$ always fits on the
right side of $R_i$, and $\tau(i)$ fits below $R_{i+1}$, so
the resulting sequences $\lambda(i) =
\attach{R_i}{\sigma(i)}{\tau(i-1)}$ are always partitions.  This
formula can be remembered by the picture 
\[ \picC{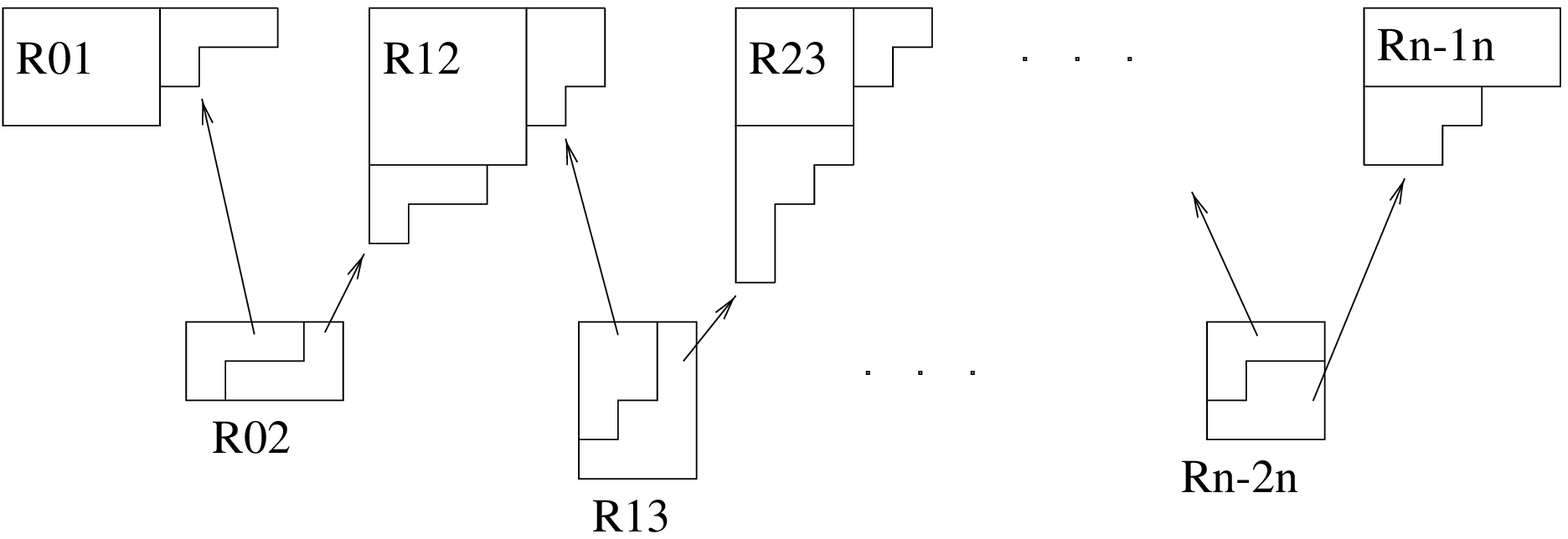} \]

The situation in the corollary covers the case of varieties of
complexes, which means that $r_{ij} = 0$ for $j-i \geq 2$.  In this
case the triangular array is
\[ \begin{matrix}
e_0 && e_1 && e_2 && e_3 && \cdots && e_n \\
& r_1 && r_2 && r_3 && \cdots && r_n \\
&& 0 && 0 && \cdots && 0
\end{matrix} \]
so the array of rectangles is
\[ \picA{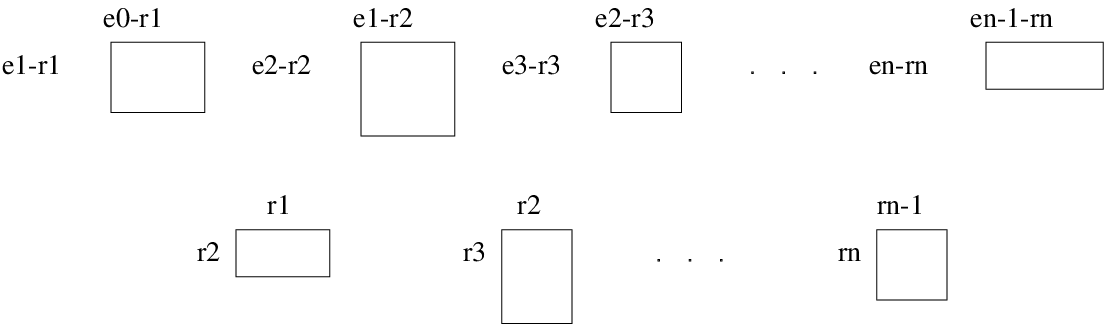} \]
P.~Pragacz reports that he had known this formula
for $\Omega_r$ in the case of varieties of complexes.

\subsection{Geometric description of $c_\lambda(r)$}
\label{sec_geomdesc}

Although we use the notation $c_\lambda(r)$ for the coefficients, it
should be emphasized that they depend only on the differences
$r_{i,j-1} - r_{ij}$ and $r_{i+1,j} - r_{ij}$, not on the integers
$r_{ij}$ themselves.  For example they are unchanged if the same
positive integer is added to each $r_{ij}$.

The linear independence of ordinary Schur polynomials $s_\lambda(x_1,
\dots, x_p)$ as $\lambda$ varies over partitions of length at most
$p$, implies that the polynomials $s_\lambda(E_\bull)$ are linearly
independent functions of the Chern classes of the bundles $E_0, E_1,
\dots, E_n$, if the ranks of the bundles are suitably large (e.g.\ if
$\ell(\lambda(i)) \leq e_i$ for $1 \leq i \leq n$).

From the preceding two paragraphs it follows that the coefficients
$c_\lambda(r)$ are uniquely determined by the geometry, i.e.\ by the
fact that $\sum c_\lambda(r) s_\lambda(E_\bull)$ is a formula for
$\Omega_r$.  To see this, one can choose the ranks $e_i$ large, and
one can find a smooth variety $X$ on which $\Omega_r$ has the expected
codimension $d(r)$, and for which the classes $s_\lambda(E_\bull)$,
for $\sum |\lambda(i)| = d(r)$, are linearly independent.  For
example, one can start with universal bundles $E_i$ on large
Grassmanians $G_i$, let $G = \prod_{i=0}^n G_i$, and set $X =
\oplus_{i=1}^n \Hom(E_{i-1}, E_i)$.

\subsection{Schubert polynomials}

In \cite{abeasis.del-fra:degenerations}, the rank conditions $r$
satisfying \refeqn{eqn_rankcond} are
described by diagrams of dots connected by lines.  One puts $e_i =
r_{ii}$ dots in column $i$, and lines are drawn between some dots in
adjacent columns.  Then $r_{ij}$ is the number of lines connecting a
dot in column $i$ to a dot in column $j$.  The example given earlier
in this section can be described by the diagram
\[ \picC{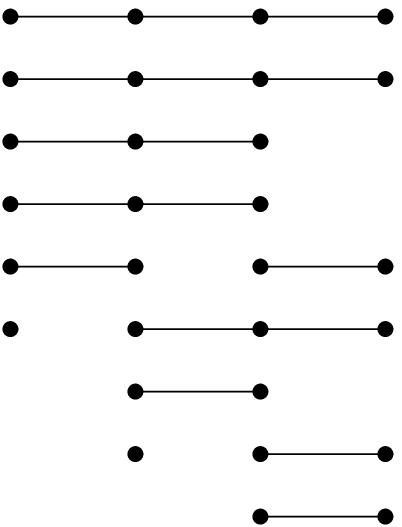} \]

Now fix a positive integer $m$.  For a permutation $w \in
S_{m+1}$, we form a diagram with $2m$ columns of lengths $1, 2, \dots,
m, m, m-1, \dots, 2, 1$.  All possible lines are drawn among the first
$m$ columns and among the last $m$ columns.  Between the two middle
columns, the $i$'th dot on the right is connected to the $w(i)$'th dot
on the left.  If $w(i) = m+1$, no connection is made.  For example, if
$m=4$ and $w = 31452$, this diagram is
\[ \picC{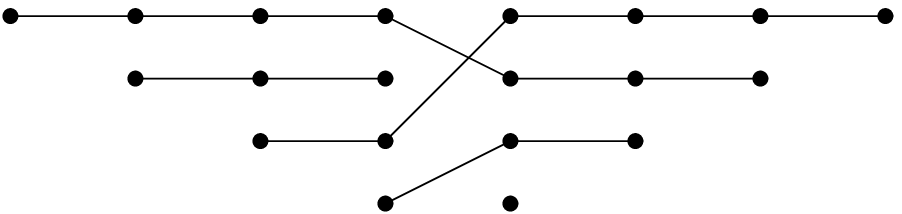} \]
The number of connections between the left column with $q$ dots and
the right column with $p$ dots is the number
\[ r_w(p, q) = \# \{ i \leq p \mid w(i) \leq q \} \,. \]
There are the maximal number of connections between two columns on the
left or between two on the right.  This means that for a sequence
$E_\bull$ of bundle maps
\[ G_1 \to G_2 \to \cdots \to G_m \to 
   F_m \to F_{m-1} \to \cdots \to F_1
\]
with $\rank(G_i) = \rank(F_i) = i$, the locus $\O_r(E_\bull)$ defined
in the introduction is exactly the
locus $\Omega_w$ defined in \cite{fulton:universal}, with the same
scheme structure.  In \cite{fulton:universal} ``universal Schubert
polynomials'' ${\mathfrak S}_w(c_\bull(F_\bull); c_\bull(G_\bull))$
were constructed, which represent the loci $\Omega_w$.  
From the fact that the formula for a locus is unique, we deduce the
following corollary.

\begin{cor}
With $r$ determined by $w$ as above, 
\[ {\mathfrak S}_w(c_\bull(F_\bull); c_\bull(G_\bull)) =
   \sum_\lambda c_\lambda(r) \, s_\lambda(E_\bull)
   \,.
\]
\end{cor}

When these bundle maps are specialized so that each $G_{i-1}
\to G_i$ is an inclusion of bundles, and each $F_i \to
F_{i-1}$ is a surjection, then ${\mathfrak S}_w(c_\bull(F_\bull);
c_\bull(G_\bull))$ becomes the double Schubert polynomial 
\[ {\mathfrak S}_w(x_1, \dots, x_m; y_1, \dots, y_m) \]
of Lascoux and Sch{\"u}tzenberger; here we set $x_i = c_1(\ker(F_i
\to F_{i-1}))$ and $y_i = c_1(G_i/G_{i-1})$.  The right side
of the formula in this corollary also simplifies in this case.  It
follows from the definition that for a partition $\tau$, we have
\[ s_\tau(G_i - G_{i-1}) = \begin{cases} 
   (y_i)^q & \text{if $\tau = (q)$, $q \geq 0$} \\
   0 & \text{otherwise}
   \end{cases}
\] 
and
\[ s_\tau(F_{i-1} - F_i) = \begin{cases}
   (-x_i)^p & \text{if $\tau = (1)^p$, $p \geq 0$} \\
   0 & \text{otherwise.}
   \end{cases}
\]
Thus $s_\lambda(E_\bull) = 0$ unless $\lambda = 
((q_2), (q_3), \dots, (q_m), \tau, (1)^{p_m}, \dots, (1)^{p_2})$, in
which case 
\[ s_\lambda(E_\bull) = 
   (-1)^{p_2 + \dots + p_m} x_2^{p_2} \cdots x_m^{p_m} y_2^{q_2}
   \cdots y_m^{q_m} s_\tau(x/y) \,,
\]
where $s_\tau(x/y) = \det(h_{\tau_i + j - i})$, $\sum h_k = \prod
(1-y_i) / \prod (1-x_j)$.  Our formula therefore writes ${\mathfrak
  S}_w(x;y)$ as a signed sum of monomials in the $x_2, \dots, x_m,
y_2, \dots, y_m$ times Schur polynomials $s_\tau(x/y)$.  When all
variables $y_i$ are set equal to zero, we have only the terms with
$q_2 = \cdots = q_m = 0$, and this writes the ordinary Schubert
polynomial ${\mathfrak S}_w(x) = {\mathfrak S}_w(x;0)$ as a signed sum
of monomials in $x_2, \dots, x_m$ times (symmetric) Schur polynomials
$s_\tau(x)$.

Unlike the inductive construction of Schubert polynomials from high
degree to low degree, our formulas are simplest for those of low
degree.

For example, for $w = 3142 \in S_4$, the corresponding array of
rectangles is
\[ \picC{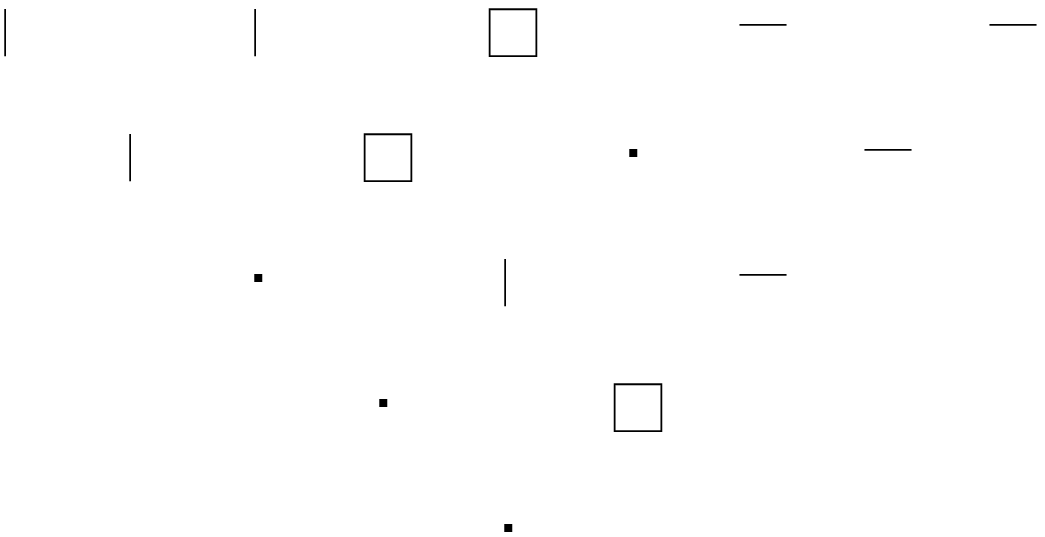} \]
Calculating $\sum c_\lambda(r) S(\lambda)$ with the algorithm of the
main theorem, working from the bottom up, one has, with $\emptyset$
the empty partition
\[ \begin{split} 
S(\emptyset) & \rightsquigarrow
S(\emptyset, 1) \\
& \rightsquigarrow
S(\emptyset, 1, \emptyset) + S(\emptyset, \emptyset, 1)  \\
& \rightsquigarrow
S(\emptyset,2,\emptyset,\emptyset) +
S(\emptyset,1,1,\emptyset) +
S(\emptyset,1,\emptyset,1) \\
& \rightsquigarrow
  {\textstyle \sum} c_\lambda(r) S(\lambda) = 
  S(\emptyset,2,1,\emptyset,\emptyset) +
  S(\emptyset,1,1\,1,\emptyset,\emptyset) + \\
&\quad\quad
  S(\emptyset,1,2,\emptyset,\emptyset) + 
  S(\emptyset,1,1,1,\emptyset) + 
  S(\emptyset,\emptyset,2\,1,\emptyset,\emptyset) + \\
&\quad\quad
  S(\emptyset,\emptyset,1\,1,1,\emptyset) +
  S(\emptyset,1,1,\emptyset,1) +
  S(\emptyset,\emptyset,1\,1,\emptyset,1) \,.
\end{split}
\]
Substituting $s_\lambda(E_\bull)$ for $S(\lambda)$, this is a formula
for the universal double Schubert polynomial associated to $w = 3142$.
It specializes to the formula
\begin{multline*}
{\mathfrak S}_{3142}(x;y) = s_{21}(x/y) + y_3 s_2(x/y) + \\
(y_3 - x_2 - x_3) s_{11}(x/y) + y_3(y_3 - x_2 - x_3) s_1(x/y)
\end{multline*}
and to 
\[ {\mathfrak S}_{3142}(x) = s_{21}(x) - (x_2 + x_3) s_{11}(x) \,. \] 

The rectangular array coming from a permutation $w \in S_{m+1}$ has
only empty rectangles and $1 \times 1$ rectangles.  In fact, this
array is determined from the diagram of the permutation denoted
$D'(w)$ in \cite[\S2]{fulton:universal}: the diagram $D'(w)$ is
reflected in a vertical line, then rotated $135^\circ$ clockwise, and
the result fitted in the bottom of the triangle; each box in $D'(w)$
is then in the position of a non-empty rectangle of the rectangular
array.  In other words the rectangle $R_{ij}$ is non-empty iff $D'(w)$
contains a box in position $(2m-j, i+1)$, which happens exactly when
$w(2m+1-j) \leq i+1$ and $w^{-1}(i+2) \leq 2m - j$.

\section{Proof of the Theorem}
\label{sec_proof}

\newcommand{\U}{\Omega^\circ}

\subsection{Geometric preliminaries}

It follows from the general discussion in \refsec{sec_quiver} that it
suffices to prove the formula for the universal locus $\O_r(E_\bull)$
in $H = \bigoplus \Hom(E_{i-1}, E_i)$.  (Throughout this section we omit
notation for pullbacks of bundles by canonical maps.)  In particular,
we know from \cite{abeasis.del-fra.ea:geometry} and
\cite{lakshmibai.magyar:degeneracy} that $\O_r(E_\bull)$ is reduced
and irreducible of the 
expected codimension $d(r)$, and that $\O_r(E_\bull)$ is the closure
of the locus $\U_r(E_\bull)$ where each of the maps $E_i(x)
\rightarrow E_j(x)$ has rank equal to $r_{ij}$ for $i < j$.  We must
prove that, with these assumptions,
\[ [\O_r(E_\bull)] = \sum_\lambda c_\lambda(r) \, s_\lambda(E_\bull) 
   \scap [H] 
\]
in the Chow group $A_{N - d(r)}(H)$, $N = \dim(H)$.

Form the Grassmannian bundle $G_0$ over $X$ with $r_i = r_{i-1,i}$ as
in \S\ref{ss_main}:
\[ G_0 = \Gr(r_1, E_1) \times_X \Gr(r_2, E_2) \times_X \cdots 
   \times_X \Gr(r_n, E_n) \,.
\]
Let $G = G_0 \times_X H$, with projection $\pi : G \rightarrow H$.
Let $0 \rightarrow A_i \rightarrow E_i \rightarrow Q_i \rightarrow 0$
be the universal exact sequences on $G_0$, and hence also on $G$.  Let
$Z \subset G$ be the intersection of the zero-schemes of the canonical
maps $E_{i-1} \rightarrow E_i \rightarrow Q_i$, i.e.\ 
\[ Z = Z(E_0 \to Q_1) \cap Z(E_1 \to Q_2) \cap
   \cdots \cap Z(E_{n-1} \to Q_n) \,.
\]
On $Z$ we have maps $E_{i-1} \to A_i$ for $1 \leq i \leq n$.
Composing these with the inclusions $A_{i-1} \subset E_{i-1}$ we get a
sequence $A_\bull$ of bundles and bundle maps on $Z$:
\[ A_1 \to A_2 \to \cdots \to A_n \,. \]
Let $\bar r$ denote the rank conditions obtained by omitting the top
row of the triangular array for $r$, and let $\O_{\bar r}(A_\bull)
\subset Z$ be the locus given by these maps and rank conditions.  It
is easy to see that $\O_{\bar r}(A_\bull)$ is mapped into
$\O_r(E_\bull)$ by $\pi$.

Now $Z$ is isomorphic to the bundle $\bigoplus_{i=1}^n \Hom(E_{i-1},
A_i)$ over $G_0$, and we have a canonical projection
\[ \rho : Z = \bigoplus_{i=1}^n \Hom(E_{i-1},A_i) \longrightarrow
   \bigoplus_{i=2}^n \Hom(A_{i-1},A_i) = H' \,.
\]
Denote by $\O'$ the universal locus $\O_{\bar r}(A_\bull)$ of $H'$.
Then $\O_{\bar r}(A_\bull)$ in $Z$ is the inverse image of $\O'$ by
$\rho$.  Since the maps on $H'$ are universal, it follows that $\O'$
is irreducible, and therefore $\O_{\bar r}(A_\bull)$ is an irreducible
subscheme of $Z$.

\[ \xymatrix{
\O' \ar@{_{(}->}[d] & \O_{\bar r}(A_\bull) \ar[l] \ar[rr] 
  \ar@{_{(}->}[d] & & \O_r(E_\bull) \ar@{_{(}->}[d] \\
H' & Z \ar[l]_\rho \ar@{^{(}->}[r] & G \ar[r]^\pi \ar[d] & H \ar[d] \\
&& G_0 \ar[r] & X
} \]

\begin{lemma}
\label{lemma_birat}
$\pi$ maps $\O_{\bar r}(A_\bull)$ birationally onto $\O_r(E_\bull)$.
\end{lemma}
\begin{proof}
Let $Z^\circ$ denote the open subset of $Z$ where the maps $E_{i-1}
\rightarrow A_i$ are surjective.  Then 
the schemes $\U_{\bar r}(A_\bull) \cap Z^\circ$ and $\U_r(E_\bull)$ are
isomorphic; they are both universal objects in the category of
schemes $Y$ over $H$, such that the pullback to $Y$ of the
sequence $E_\bull$ satisfies $\rank(E_i \to E_j) = r_{ij}$ everywhere
on $Y$.  Since both $\O_{\bar r}(A_\bull)$ and $\O_r(E_\bull)$ are
irreducible schemes, and since $\U_r(E_\bull) \neq \emptyset$, the
assertion follows.
\end{proof}

Note that this argument gives a direct proof that the codimension of
$\O_r(E_\bull)$ in $H$ is $d(r)$.  Indeed, by induction we know that
$\O'$ has codimension $d(\bar r)$ in $H'$, and $\Omega_{\bar r}
(A_\bull)$ must then have the same codimension in $Z$.  We conclude
that the codimension of $\Omega_r(E_\bull)$ in $H$ is
\[ d(\bar r) + \dim(H) - \dim(Z) = d(r) \,. \]

\subsection{Proof of the main theorem}
By induction on $n$ we know that $[\O'] = \sum c_\mu(\bar r)
s_\mu(A_\bull) \scap [H']$, so 
\[ [\O_{\bar r}(A_\bull)] = \rho^*[\O'] = 
   \sum_\mu c_\mu(\bar r) \, s_\mu(A_\bull) \scap [Z] \,.
\]
Furthermore $[Z] = \prod_{i=1}^n c_{\text{top}}( \Hom(E_{i-1}, Q_i) )
\scap [G]$, so
\[ [Z] = \prod_{i=1}^n s_{R'_i}(Q_i - E_{i-1}) \scap [G] \,, \]
where $R'_i = (e_{i-1})^{e_i - r_i}$.  Since $\pi_*[\O_{\bar
  r}(A_\bull)] = [\O_r(E_\bull)]$ by \reflemma{lemma_birat}, we are
therefore reduced to proving the identity
\begin{equation}
\label{eqn_toprove}
\pi_*\left( \sum_\mu c_\mu(\bar r) s_\mu(A_\bull) \cdot 
\prod_{i=1}^n s_{R'_i}(Q_i - E_{i-1}) \scap [G] \right)
= \sum_\lambda c_\lambda(r) s_\lambda(E_\bull) \scap [H] \,.
\end{equation}
For this we need the following basic Gysin formula of Pragacz
\cite[Prop.~2.2]{pragacz:enumerative}, whose proof comes from \cite{jozefiak.lascoux.ea:classes}, cf.\ \cite[App.~F]{fulton.pragacz:schubert}. 
\begin{lemma}
\label{lemma_gysin}
Let $E$ and $F$ be vector bundles of ranks $e$ and $f$ on a variety
$X$.  Let $0 \leq d \leq \min(e,f)$.  Let $G = \Gr(d, F)$ be the
Grassmann bundle, with projection $\pi: G \rightarrow X$ and universal
exact sequence $0 \rightarrow A \rightarrow F \rightarrow Q
\rightarrow 0$.  Let $q = f - d$, $R = (e-d)^q$, and $R' = (e)^q$.
For any partitions $\lambda$ and $\mu$, with $\lambda$ of length at
most $q$, and any $\alpha \in A_*(X)$, 
\[ \pi_*(s_{R' + \lambda}(Q - E) s_\mu(A - E) \scap \pi^* \alpha)
   = s_{R + \lambda, \mu}(F - E) \scap \alpha \,.
\]
\end{lemma}
We also need the following special case of the factorization formula
of Lascoux and Sch{\"u}tzenberger
\cite{lascoux.schutzenberger:polynomes} and Berele and Reger
\cite{berele.regev:hook}, cf.\ \cite{pragacz:enumerative}.
\begin{lemma}
\label{lemma_factor}
Let $E$ and $F$ be vector bundles of ranks $e$ and $f$.  Let $R =
(e)^f$.  Let $\lambda$ be a partition of length at most $f$.  Then
\[ s_\lambda(F) s_R(F - E) = s_{R + \lambda}(F - E) \,. \]
\end{lemma}
Note that this identity follows from \reflemma{lemma_gysin}.  Finally
we need the basic identity \cite[\S1.5]{macdonald:symmetric*2}:
\begin{lemma}
\label{lemma_split}
For bundles $E_1$, $E_2$, and $E_3$, and a partition $\mu$, 
\[ s_\mu(E_3 - E_1) = 
   \sum c^\mu_{\sigma\tau} s_\sigma(E_2 - E_1) s_\tau(E_3 - E_2) \,,
\]
the sum over partitions $\sigma$ and $\tau$ with $|\sigma| + |\tau| =
|\mu|$, with $c^\mu_{\sigma \tau}$ the Littlewood-Richardson
coefficient.
\end{lemma}

Now we can prove \refeqn{eqn_toprove}.  First use
\reflemma{lemma_split} to replace each factor $s_{\mu(i)}(A_{i+1} -
A_i)$ that occurs in each $s_\mu(A_\bull)$ on the left side of
\refeqn{eqn_toprove} by the sum
\begin{multline*}
\sum c^{\mu(i)}_{\sigma(i),\tau(i)} \, s_{\sigma(i)}(E_i - A_i) 
s_{\tau(i)}(A_{i+1} - E_i) \\
= \sum c^{\mu(i)}_{\sigma(i),\tau(i)} \, s_{\sigma(i)}(Q_i)
s_{\tau(i)}(A_{i+1} - E_i) \,.
\end{multline*}
Note that $s_{\sigma(i)}(Q_i) = 0$ if $\ell(\sigma(i)) > \rank(Q_i) =
e_i - r_i$.  Next use \reflemma{lemma_factor} to replace each
$s_{\sigma(i)}(Q_i) \cdot s_{R'_i}(Q_i - E_{i-1})$ in the result by
$s_{R'_i+\sigma(i)}(Q_i - E_{i-1})$.  The left side of
\refeqn{eqn_toprove} becomes
\begin{multline*}
 \sum_\mu c_\mu(\bar r) \sum_{\sigma(i),\tau(i)} 
   \left( \prod_{i=1}^{n-1} c^{\mu(i)}_{\sigma(i),\tau(i)} \right)
   \cdot \\
   \pi_*\left( \prod_{i=1}^n s_{R'_i+\sigma(i)}(Q_i - E_{i-1})
   s_{\tau(i-1)}(A_i - E_{i-1})  \scap [G] \right) \,.
\end{multline*}
Finally, $n$ applications of \reflemma{lemma_gysin} yields
\begin{multline*}
\pi_*\left( \prod_{i=1}^n s_{R'_i+\sigma(i)}(Q_i - E_{i-1})
s_{\tau(i-1)}(A_i - E_{i-1})  \scap [G] \right) \\
= \prod_{i=1}^n s_{R_i+\sigma(i),\tau(i-1)}(E_i - E_{i-1}) \scap [H] \,,
\end{multline*}
and this gives the required formula
\begin{multline*}
\sum_\mu c_\mu(\bar r) \sum_{\sigma(i), \tau(i)}
\left( \prod_{i=1}^{n-1} c^{\mu(i)}_{\sigma(i),\tau(i)} \right)
\prod_{i=1}^n s_{R_i+\sigma(i),\tau(i-1)}(E_i - E_{i-1}) \scap [H] \\
= \sum_\lambda c_\lambda(r) s_\lambda(E_\bull) \scap [H] \,.
\end{multline*}

Although we have stated it for varieties over a field, the theorem
(and its proof) extend readily to schemes of finite type over a
regular base, as in \cite[\S20]{fulton:intersection}.

\subsection{A generalization}
\label{sec_paths}

There is a generalization of the theorem which may be useful in its
own right, and which gives some insight into the proof.  (It is not
needed in this paper.)
Fix $E_0, \dots, E_n$ on $X$, and $r = (r_{ij})$ satisfying
\refeqn{eqn_rankcond}.  Let $H = \bigoplus \Hom(E_{i-1}, E_i)$, on
which the tautological bundle maps are universal, and one has the
universal locus $\O_r \subset H$.

Let $\pi : F \to H$ be the partial flag bundle parameterizing flags in
each $E_j$ of ranks $r_{0j}, r_{1j}, \dots, r_{j-1,j}$.  Let $E_{0j}
\subset E_{1j} \subset \cdots \subset E_{j-1,j} \subset E_j$ denote
the tautological flags of vector bundles on $F$, $1 \leq j \leq n$.

Let $Z \subset F$ be the locus on which the image of $E_{i,j-1}$ by
the map $E_{j-1} \to E_j$ is contained in the subbundle $E_{ij}$,
i.e.\ $Z$ is the subscheme defined by the vanishing of all maps
$E_{i,j-1} \to E_j/E_{ij}$ for $i < j$.

One sees as in \reflemma{lemma_birat} that $\pi$ maps $Z$ birationally
onto $\O_r$.
In fact, if $X$ is non-singular, this construction gives a canonical
resolution of singularities of the universal locus $\O_r$.  It is 
easy to see that the class of $Z$ is given by
\[ [Z] = \prod_{i<j} z_{ij} \,, \]
where $z_{ij} = c_{\text{top}}(\Hom(E_{i,j-1}, E_{i+1,j}/E_{ij}))$.

Consider a path $\gamma$ through the triangular array for $r$, going
from $r_{00}$ to $r_{nn}$.  The path must be a union of line segments
between neighboring rank conditions, and it must intersect any
vertical line at most once.
\[ \picC{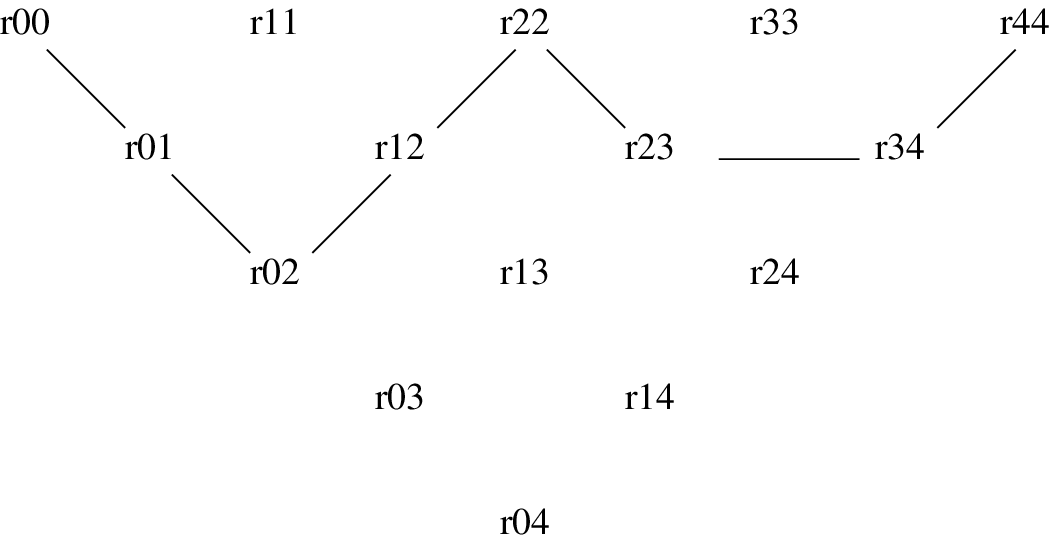} \]
For each $j = 1, \dots, n$, let $k_j$ be minimal such that $\gamma$ goes
through $r_{k_j,j}$.  
Let $F(\gamma)$ be the partial flag bundle parameterizing flags in
$E_j$ of ranks $r_{k_j,j}, r_{k_j+1,j}, \dots, r_{j-1,j}$, and let
$E_{k_j,j} \subset  \cdots \subset E_{j-1,j} \subset E_j$ be the
tautological flags on $F(\gamma)$, $1 \leq j \leq n$.  If the path has
$m$ line segments, we let $A_0, A_1, \dots, A_m$ denote the vector
bundles on $F(\gamma)$ corresponding to the rank conditions passed
through by the path.  (In the illustration, $m = 7$, and the bundle
sequence is $E_{00}, E_{01}, E_{02}, E_{12}, E_{22}, E_{23}, E_{34},
E_{44}$.) 

Let $\O_r(\gamma) \subset F(\gamma)$ be the subscheme defined by the
conditions that each map $E_{i,j-1} \to E_j/E_{ij}$ vanishes for
$r_{ij}$ on or above the path, and $\rank(E_{ip} \to E_j) \leq r_{ij}$
for $r_{ip}$ on or above the path and $p \leq j$.  The canonical maps $F \to
F(\gamma) \to X$ map $Z$ birationally onto $\O_r(\gamma)$ which in
turn is mapped birationally onto $\O_r$.  Our goal is to give a
formula for the class of $\O_r(\gamma)$.  To do this, we define a
formal linear combination $\Phi(\gamma)$ of symbols $S(\lambda)$,
where $\lambda$ is a sequence of partitions, one for each line segment
in $\gamma$.  The formula for $\O_r(\gamma)$ is obtained by replacing
each $S(\lambda(1), \dots, \lambda(m))$ by $\prod_{i=1}^m
s_{\lambda(i)}(A_i - A_{i-1})$, and multiplying the result by $\prod
z_{ij}$, the product over all $i,j$ such that $r_{ij}$ is on or above
$\gamma$.

We define $\Phi(\gamma)$ inductively.  If $\gamma$ is the lowest possible path,
going from $r_{00}$ to $r_{0n}$ to $r_{nn}$, then $\Phi(\gamma) =
S(\emptyset, \dots, \emptyset)$, where the empty partition $\emptyset$
is repeated $2n$ times.  Otherwise, we can find a path $\gamma'$ that is
equal to $\gamma$ except it goes lower at one place, in one of the
following ways:
\begin{center}
\begin{tabular}{ccccc}
& \hspace{0.5cm} & \hspace{0.35cm}$\gamma'$ & \hspace{1cm} 
& \hspace{0.35cm}$\gamma$
\vspace{0.65cm} \\
Case 1: 
\vspace{-0.75cm} \\
&& \picC{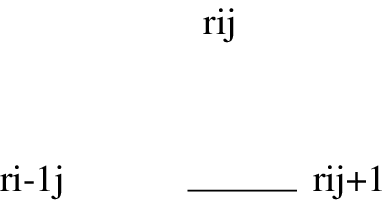} && \picC{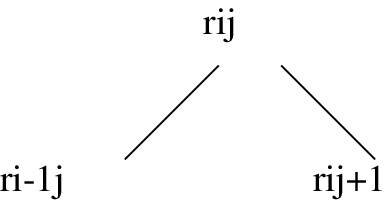}
\vspace{0.65cm} \\
Case 2:
\vspace{-0.75cm} \\
&& \picC{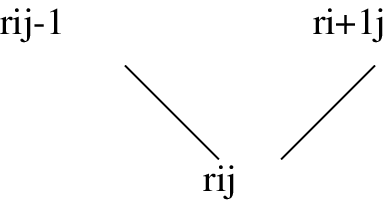} && \picC{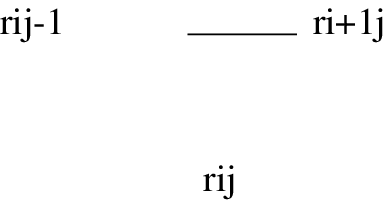}
\end{tabular}
\vspace{0.2cm}
\end{center}

In Case 1, we obtain $\Phi(\gamma)$ from $\Phi(\gamma')$ by replacing
each symbol $S(\dots, \mu, \dots)$ by $\sum c^\mu_{\sigma,\tau}
S(\dots, \sigma, \tau, \dots)$.  Note that in this case we have
$F(\gamma) = F(\gamma')$ and $\O(\gamma) = \O(\gamma')$.  For Case 2,
each symbol $S(\dots, \mu, \nu, \dots)$ in $\Phi(\gamma')$ is replaced
by the symbol $S(\dots, \attach{R_{ij}}{\nu}{\mu}, \dots)$.  Note that
the partitions $\lambda(i)$ are always empty for line segments on the
left or right edges of the triangle.

The proof that this polynomial gives a formula for $[\O_r(\gamma)]$ is
similar to that of our main theorem: one shows that $[Z]$ pushes
forward to the class of this polynomial.  The induction step, in
either Case 1 or 2, is more transparent, as changes are made in only
one segment of the formula.

If bundle maps and flags of subbundles are given on $X$ corresponding
to ranks on or above the path $\gamma$, these determine a section $s :
X \to F(\gamma)$, so one has corresponding formulas for the classes
$s^*[\O_r(\gamma)]$.  When $\gamma$ is the horizontal path across the
top of the diagram, we recover the main theorem.

\section{On the coefficients $c_\lambda(r)$}
\label{sec_coefs}

There are some properties of the coefficients $c_\lambda(r)$ that
follow from geometry, i.e.\ from the main theorem, although they are
not obvious from the algorithm defining them.  We describe these
first, and then discuss properties we believe for other reasons.
We conclude with a comparison of the numbers $c_\lambda(r)$ with
Littlewood-Richardson coefficients.

Consider the dual sequence
\[ E_n^\vee \rightarrow E_{n-1}^\vee \rightarrow \cdots 
   \rightarrow E_0^\vee 
\]
with dual rank conditions, $\rank(E_j^\vee \rightarrow E_i^\vee) \leq
r_{ij}$, which we denote by $r^\vee$; then
$\Omega_{r^\vee}(E_\bull^\vee) = \Omega_r(E_\bull)$.  Note that the
rectangular array for $r^\vee$ is obtained by reflecting that for $r$
in a vertical line, and replacing each rectangle by its transpose.
Using the basic identity that $s_\lambda(F-E) = s_{\lambda'}(E^\vee -
F^\vee)$, where $\lambda'$ is the transpose of $\lambda$, we find that
\begin{equation}
\label{eqn_duality}
c_{\lambda^\vee}(r^\vee) = c_\lambda(r) \,,
\end{equation}
where, if $\lambda = (\lambda(1), \dots, \lambda(n))$, we put
$\lambda^\vee = (\lambda(n)', \dots, \lambda(1)')$.

It can happen that for some $k$, all of the rank conditions $\rank(E_i
\rightarrow E_k) \leq r_{ik}$ and $\rank(E_k \rightarrow E_j) \leq
r_{kj}$ follow from other rank conditions.  This happens when, in the
rectangle diagram, all the rectangles on the two $45^\circ$ lines
descending from position $k$ are empty.  For the example $G_1
\rightarrow G_2 \rightarrow G_3 \rightarrow F_3 \rightarrow F_2
\rightarrow F_1$ considered at the end of \refsec{sec_quiver}, with
rank conditions 
$r$ coming from $w = 3142 \in S_4$, the bundles $G_1$ and $F_2$ are
inessential in this way.  If an inessential bundle $E_k$ is omitted,
one has a shorter sequence $E'_\bull : E_0 \to \cdots \to E_{k-1} \to
E_{k+1} \to \cdots \to E_n$, with the map $E_{k-1} \rightarrow E_{k+1}$
being $\phi_{k+1} \circ \phi_k$, and corresponding rank conditions
$r'$; the array of rectangles for $r'$ is obtained by 
omitting the $45^\circ$ lines of empty rectangles and moving all
rectangles below up a row.  For example, if $G_1$ and $F_2$ are
omitted from the example, one has $G_2 \rightarrow G_3 \rightarrow F_3
\rightarrow F_1$, with rectangular array
\[ \picC{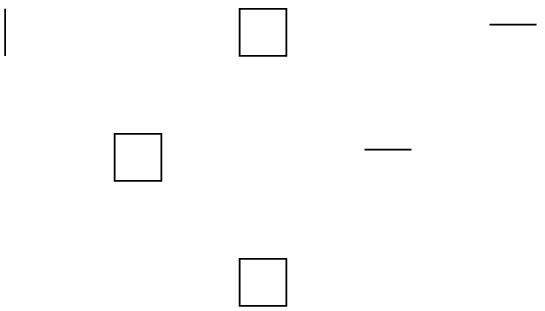} \]
\reflemma{lemma_split} can be used to expand any $s_\rho(E_{k+1} -
E_{k-1})$ occurring in the formula $\sum c_\mu(r') s_\mu(E'_\bull)$ as
a sum of terms of the form $s_\sigma(E_k - E_{k-1}) \cdot
s_\tau(E_{k+1} - E_k)$.  Since $\Omega_{r'}(E'_\bull) =
\Omega_r(E_\bull)$, with this interpretation we have
\begin{equation}
\label{eqn_omit}
\sum_\mu c_\mu(r') s_\mu(E'_\bull) = 
\sum_\lambda c_\lambda(r) s_\lambda(E_\bull) \,.
\end{equation}

Now we turn to our conjectured formula for the coefficients
$c_\lambda(r)$, which interprets them by counting Young tableaux, in a
way similar to and generalizing the classical Littlewood-Richardson
rule.  Recall that a (semistandard) {\em Young tableau\/} is a filling
of the boxes in the Young diagram of a partition with integers that
are weakly increasing in rows and strictly increasing down columns.
Two Young tableaux $P$ and $Q$ can be multiplied to give another Young
tableau denoted $P \cdot Q$.  One way to do this is to arrange $P$ and
$Q$ corner to corner and play the jeu de taquin, sliding inside
corners through but keeping the weak and strict orderings.  For
example, if $P = \picC{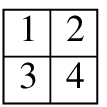}$ and $Q =
\picC{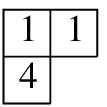}$, this can be carried out by the sequence of
moves
\[ \picC{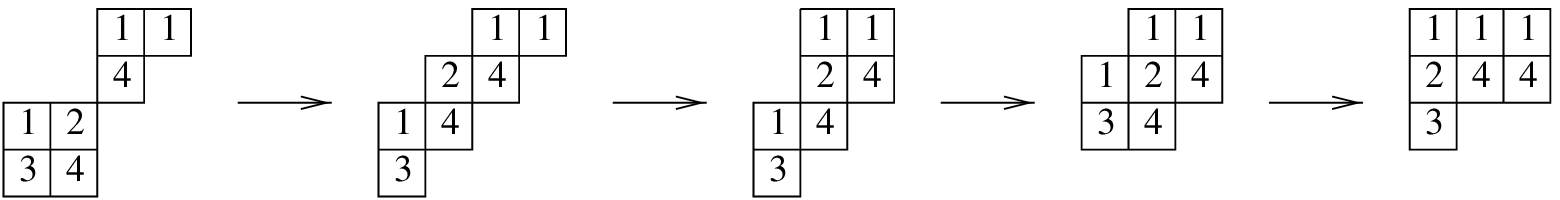} \]
The final tableau is $P \cdot Q$.  The main fact is that this product
is independent of choice of the sequence of inside corners, from which
it follows that the set of tableaux form an associative monoid, called
the plactic monoid.  With this notion, the Littlewood-Richardson
number $c^\lambda_{\mu,\nu}$ is the number of ways a given tableau $T$
of shape $\lambda$ can be factored into a product $T = P \cdot Q$,
such that $P$ has shape $\mu$ and $Q$ has shape $\nu$.  For proofs and
relations with Schur polynomials, see
\cite{lascoux.schutzenberger:mono} or \cite{fulton:young}.

Given rank conditions $r$ (satisfying \refeqn{eqn_rankcond} as
always), form the array of rectangles $R_{ij}$.  We choose a fixed
tableau $T_{ij}$ on each shape $R_{ij}$, with the condition that each
entry of $T_{ij}$ must be strictly smaller than any entry of $T_{kl}$
if $R_{kl}$ lies in the wedge cut out by $45^\circ$ lines below
$R_{ij}$, i.e.\ if $k \leq i$ and $l \geq j$ with $(k,l) \neq (i,j)$.

From this array of rectangular tableaux we will construct a set of
$n$-tuples of tableaux $(T_1, \dots, T_n)$ that we call {\em factor
  sequences}.  Our conjecture is that $c_\lambda(r)$ is the number of
factor sequences $(T_1, \dots, T_n)$ such that $T_i$ has shape
$\lambda(i)$ for $1 \leq i \leq n$.  We first explain this for $n =
3$, where we start with an array of rectangular tableaux:
\[ \begin{matrix}
A && B && C \\
& D && E \\
&& F
\end{matrix} \]
Factor $F$ into a product $F = F_1 \cdot F_2$ of tableaux.  Pass $F_1$
up to the left, and multiply it to $D$ from the right.  Pass $F_2$ up
to the right and multiply it to $E$ from the left.  Then factor the
results:
\begin{align*}
D \cdot F_1 &= D_1 \cdot D_2 & &\text{and}& 
F_2 \cdot E &= E_1 \cdot E_2 \,. 
\end{align*}
Pass the results up to the left and right, arriving at tableaux $A
\cdot D_1$, $D_2 \cdot B \cdot E_1$, and $E_2 \cdot C$.  This gives a
factor sequence $(T_1, T_2, T_3) = (A \cdot D_1,\, D_2 \cdot B \cdot
E_1,\, E_2 \cdot C)$.

In general one proceeds by induction.  A factor sequence for the given
array of rectangular tableaux is obtained by forming a factor sequence
$(S_1, \dots, S_{n-1})$ for the array of the bottom $n-1$ rows.
Factor each $S_i$ arbitrarily into $S_i = P_i \cdot Q_i$.  Then
\[
(T_1,\dots,T_n) = (T_{01} \cdot P_1, \dots, Q_{i-1} \cdot
T_{i-1,i} \cdot P_i, \dots, Q_{n-1} \cdot T_{n-1,n})
\]
is a factor sequence for the given array.

\begin{conj}
$c_\lambda(r)$ is the number of factor sequences $(T_1, \dots, T_n)$
of shape $\lambda = (\lambda(1), \dots, \lambda(n))$ that can be made
from a given array of rectangular tableaux.
\end{conj}

The conjecture has a number of consequences:
\begin{C1}
Each $c_\lambda(r)$ is a non-negative integer.
\end{C1}
\begin{C2}
The coefficients $c_\lambda(r)$ depend only on the rectangles $R_{ij}$, not on their sides.
\end{C2}
This means that if one of the sides of a rectangle $R_{ij}$ is $0$,
the length of the other side is irrelevant.  (The algorithm of the
main theorem shows this when the height of a rectangle is 0, but not
when the width is 0.)

Implicit in the conjecture is the assertion
\begin{C3}
\label{cons_indep}
The number of factor sequences of shape $\lambda$ is independent
  of choice of fixed tableaux $T_{ij}$.
\end{C3}

Granting C3, it is not hard to see that the
conjectured formula for the $c_\lambda(r)$ satisfies the duality
\refeqn{eqn_duality}.  For this one chooses the $T_{ij}$ so that no
entry appears more than once, and uses the fact that factoring a
tableau $T$ with distinct entries into $P \cdot Q$ is equivalent to
factoring its conjugate $T'$ into $Q' \cdot P'$.  It is also not
difficult to verify that the conjectured formula satisfies the
property \refeqn{eqn_omit} for omitting inessential bundles.

The conjecture is true for the case where $R_{ij}$ is empty for
$j-i>2$.  This follows from the description in \refcor{cor_complexes}
of \refsec{sec_quiver}, together with the fact that for a tableau $T$
of rectangular 
shape $R$, for each $\sigma$ and $\tau$ that fit together to form $R$,
there is a unique factorization $T = P \cdot Q$ with $P$ of shape
$\sigma$ and $Q$ of shape $\tau$; conversely, any factorization of $T$
has factors of shapes that fit together to form $R$.

The conjecture has been proved when $n \leq 3$.  More generally, it is
proved when $R_{ij}$ is empty for $j - i > 3$ and no two non-empty
rectangles in the third row are adjacent.  The proof depends on a
wonderful involution on pairs of tableaux produced for us by S.~Fomin.
This proof is given in \cite{buch:combinatorics}.

One reason that the combinatorial formula is hard to work with is that
a given factor sequence can arise in many ways by the algorithm that
produces them.  At first glance it would appear that to tell if some
$(T_1, \dots, T_n)$ is a factor sequence, one would have to test all
possible ways of carrying out the sequence of factorings.  However,
there is a direct test.  For this, define $P_i$ to be the part of
$T_i$ lying to the right of the rectangle $R_i$, and define
$Q_{i-1}$ to be everything lying below $T_{i-1,i}$ and $P_i$:
\[ \picC{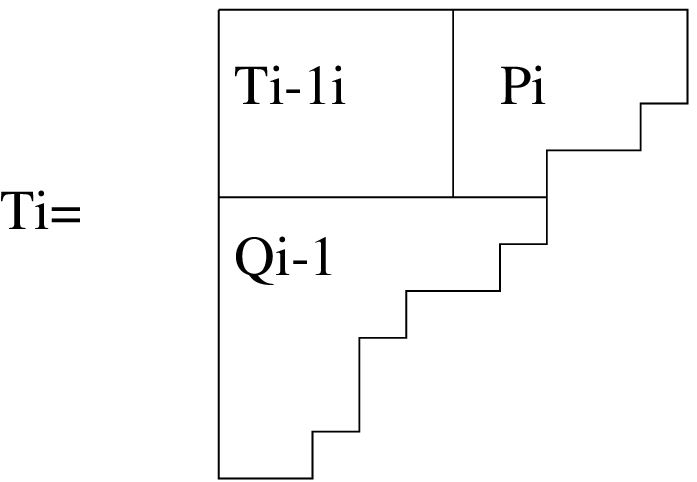} \]
Then $(T_1, \dots, T_n)$ is a factor sequence if and only if $Q_0$ and
$P_n$ are empty and $(P_1 \cdot Q_1, \dots, P_{n-1} \cdot Q_{n-1})$ is
a factor sequence for the lower $n-1$ rows of the array.  By induction
this gives a direct algorithm to test, from the top down.  Note that
this algorithm, like the theorem, uses the height of a rectangle
$R_{ij}$ even if its width is zero.  This criterion is proved in
\cite{buch:combinatorics}.

The full conjecture follows from an assertion that Fomin's involution
preserves factor sequences.  This assertion is true for $n \leq 3$,
and it has been verified in 500,000 randomly generated examples for
$n \leq 10$.  For a discussion of Fomin's involution and the
discussion of this, we refer again to \cite{buch:combinatorics}.

The numbers $c_\lambda(r)$ generalize Littlewood-Richardson numbers
$c^\gamma_{\alpha,\beta}$ in fact as well as in spirit.  To see this,
take any rectangle $R$ containing $\gamma$, and let $\tau$ be the
complement of $\gamma$ in $R$ (rotated $180^\circ$):
\[ \picC{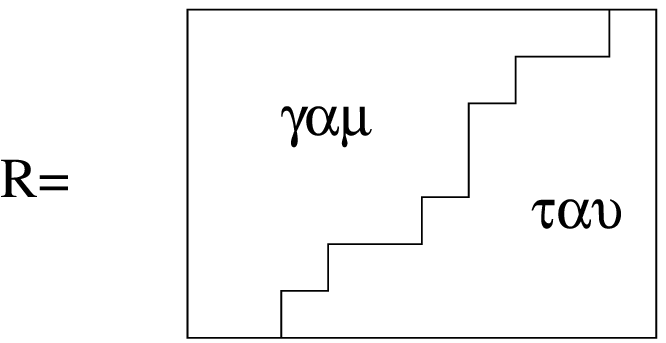} \]
Form the array of rectangles
\[ \picC{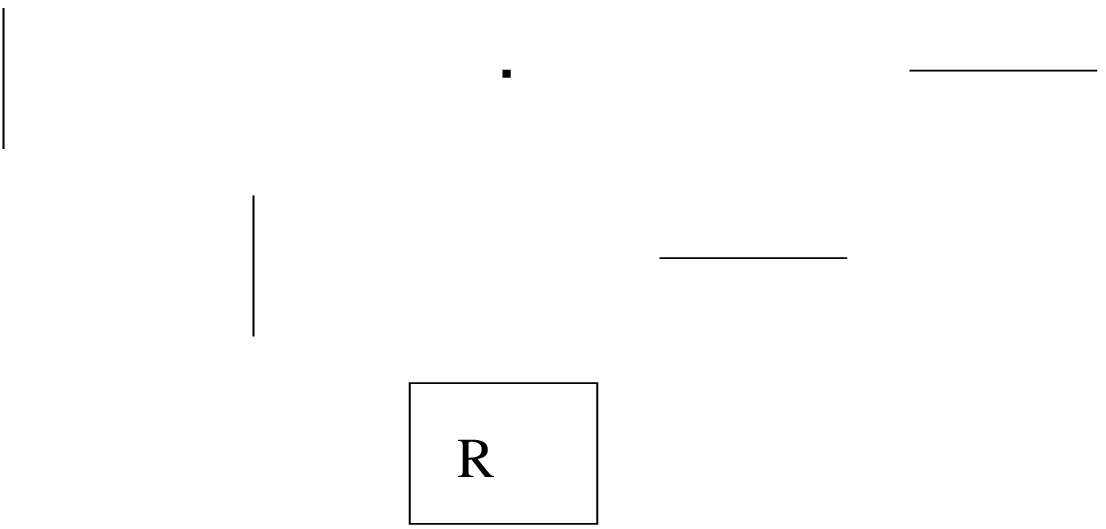} \]
where the vertical lines are the height of $R$, the horizontal lines
are its width, and the dot is empty.  Choose $r$ giving rise to this
array.  Then
\[ c^\gamma_{\alpha,\beta} = c_\lambda(r) \,, \] 
with $\lambda = (\alpha, \beta, \tau)$.  This follows easily from the
theorem (and also from the conjecture).

There are analogous conjectures for the coefficients of the
polynomials for the more general loci described in \S\ref{sec_paths}.
In particular, all these coefficients should also be positive.
Details will be given in \cite{buch:combinatorics}.

\providecommand{\bysame}{\leavevmode\hbox to3em{\hrulefill}\thinspace}

\end{document}